\documentclass{article}%
\usepackage{amsmath}%
\usepackage{amsfonts}%
\usepackage{amssymb}%
\usepackage{graphicx}

\begin{document}

\title{Recursion Relations and Functional Equations for the Riemann Zeta Function}
\author{Henrik Stenlund\thanks{The author is grateful to Visilab Signal Technologies for supporting this work.}
\\Visilab Signal Technologies Oy, Finland}
\date{July 17th, 2011}
\maketitle

\begin{abstract}
New recursion relations for the Riemann zeta function are introduced. Their derivation started from the standard functional equation. The new functional equations have both real and imaginary increment versions and can be applied over the whole complex plane. We have developed various versions of the recursion relations eliminating each of the coefficient functions, leaving plain zeta functions.
\footnote{Visilab Report \#2011-07. Revision 3,  08-09-2011, revised with two new recursion relations to section 4.}
\subsection{Keywords}
Riemann zeta function, zeros of zeta function, recursion relation of zeta function, functional equation of zeta function
\subsection{Mathematical Classification}
Mathematics Subject Classification 2010: 11M26, 11E45, 11M41, 11S40, 30G30
\end{abstract}

\section{Introduction}
\subsection{Known Recursion Formulas}
\noindent The functional equation of the Riemann zeta function with the argument $s \in C$ has been derived by Riemann \cite{Riemann1858}
\begin{equation}
\zeta(s)=2^{s}\pi{^{s-1}}sin(\frac{\pi{s}}{2})\Gamma(1-s)\zeta(1-s) \ \  s \in C \label{eqn10}
\end{equation}
This equation is symmetric in change of variable $s\rightarrow1-s$ as is easy to prove. There are only a very few recursion relations known for the zeta function. They are mostly simple restatements of equation (\ref{eqn10}) (textbooks \cite{Patterson1999}, \cite{Edwards2001} and \cite{Titchmarsh1999}). The following equation (\ref{eqn11}) is one of them and can be reached by elementary operations.
\begin{equation}
\zeta(1-s)=\frac{2}{(2\pi)^s}cos(\frac{\pi{s}}{2})\Gamma(s)\zeta(s)  \label{eqn11}
\end{equation}
The third non-trivial functional equation is equation (\ref{eqn17}) below. These are the only known recursion relations without excessive complications, like integration or summation. Recursion relations with summations are handled e.g. in \cite{Bi2010}, \cite{Apostol1954}, \cite{Baez-Duarte2010}, \cite{Adamchik1998}, \cite{Choi1995},  \cite{Rubinstein2009} and \cite{Liu2004} and thus outside the scope of this work. Some new ideas have appeared too. Lagarias \cite{Lagarias2006} studied the Riemann $\xi(s)$, developed originally by Riemann.  
\begin{equation}
\xi(s)=\frac{1}{2}s(s-1)\pi^{(\frac{-s}{2})}\Gamma(\frac{s}{2})\zeta(s)  \label{eqn12}
\end{equation}
It has a functional equation
\begin{equation}
\xi(s)=\xi(1-s)  \label{eqn14}
\end{equation}
with an obvious symmetry. Lagarias developed some recursions for this function which is a bit different from the zeta function in behavior. Ossicini \cite{Ossicini2010} has treated similar functions.

Tyagi \cite{Tyagi2007} derives a new integral formula for the zeta function, a new recursive summation formula and also displays a simple functional equation
\begin{equation}
\Gamma(\frac{s}{2}-1)\pi^{(\frac{-s}{2})}\zeta(s)=\Gamma(\frac{1-s}{2}-1)\pi^{\frac{-(1-s)}{2}}\zeta(1-s) \ (FALSE)  \label{eqn16}
\end{equation}
Unfortunately, this equation is false. The correct form should be
\begin{equation}
\Gamma(\frac{s}{2})\zeta(s)=\Gamma(\frac{1-s}{2})\pi^{s-\frac{1}{2}}\zeta(1-s)  \label{eqn17}
\end{equation}
which is a known variant of the basic functional equation (\ref{eqn10}). It can be derived easily by applying the Legendre duplication formula. 

A recent article \cite{Chiang2006} proves that certain types of recursion relations are not possible. The relations we are to develop do not fall into that category. Baez-Duarte has attempted to systematize the functional equations \cite{Baez-Duarte2002}.

\subsection{Developing New Recursion Relations}
A recursion relation means an algebraic relationship between values of a function at various points over the argument space. For complex variables, this would mean a group of points over the whole complex plane. In general, the terms may have any multiplicative factors, either constant or any other function. In a more general case, the arguments may have any scaling, including reflection. Extending generalization further leads to combinations of functions with any argument values and with any coefficient function. The functional equations (\ref{eqn10}), (\ref{eqn11}) and (\ref{eqn17}) are examples having these features.

The motivation for finding new recursion relations is two-fold. First, the original equation (\ref{eqn10}) is complicated in behavior caused by the $2^s$ and $sin(\frac{\pi{s}}{2})$ functions together with the $\Gamma(1-s)$ function. They make it very obscure to see the function's behavior when the argument $s$ varies over the complex plane. Especially studying nontrivial zeros of the zeta function is difficult. Elimination of these functions from the functional equation may be useful. The second problem is the entire missing of imaginary increment recursions. 

In the following we derive recursion relations for the Riemann zeta function from equation (\ref{eqn10}). They will be equivalent to the original functional equation. The zeta function is known to be analytic everywhere and has a single pole at $s=1$. We perform only simple operations preserving analyticity throughout to the final functional equations. This applies to all operations in this work and all equations are analytic unless indicated otherwise. In the next section we present the resulting functional equations eliminating all unwanted functions. The imaginary-valued argument increments are treated in a similar way in the following section. In the last section we show some additional recursion relations and results. 

\section{Functional Equations, Real Increments}
Starting from equation (\ref{eqn10}) with an increment to eliminate the $\Gamma(1-s)$ function, in a simple manner, gives the following.
\begin{equation}
\frac{\zeta(s)}{\zeta(1-s)}=-s\frac{tan(\frac{\pi{s}}{2})\zeta(1+s)}{2\pi\zeta(-s)} \  \   s \in C  \label{eqn20}
\end{equation}
We managed to eliminate the $2^s$ function at the same time. We continue by eliminating the $tan()$ function and get the following 
\begin{equation}
\zeta(s)=-s(s+1)\frac{\zeta(1-s)\zeta(s+2)}{4\pi^2\zeta(-1-s)}  \label{eqn30}
\end{equation}
We still have a polynomial in $s$ and need to get rid of it yielding
\begin{equation}
2s^2=-\frac{\zeta(s)\zeta(-1-s)4\pi^2}{\zeta(1-s)\zeta(s+2)}-\frac{\zeta(s-1)\zeta(-s)4\pi^2}{\zeta(2-s)\zeta(s+1)}  \label{eqn40}
\end{equation}
Now there is a square of $s$ remaining. The next step leads to the equation
\begin{equation}
2s=-\frac{\zeta(s)\zeta(-1-s)4\pi^2}{\zeta(1-s)\zeta(s+2)}+\frac{\zeta(s-1)\zeta(-s)4\pi^2}{\zeta(2-s)\zeta(s+1)}  \label{eqn50}
\end{equation}
The remaining $s$ is eliminated and we will have
\begin{equation}
\frac{1}{2\pi^2}=2\frac{\zeta(s)\zeta(-1-s)}{\zeta(1-s)\zeta(s+2)}-\frac{\zeta(s+1)\zeta(-2-s)}{\zeta(-s)\zeta(s+3)}-\frac{\zeta(s-1)\zeta(-s)}{\zeta(2-s)\zeta(s+1)}    \label{eqn60}
\end{equation}
The constant term is the last one and we get
\begin{equation}
0=-3\frac{\zeta(s)\zeta(-1-s)}{\zeta(1-s)\zeta(s+2)}+3\frac{\zeta(s+1)\zeta(-2-s)}{\zeta(-s)\zeta(s+3)}+\frac{\zeta(s-1)\zeta(-s)}{\zeta(2-s)\zeta(s+1)}-\frac{\zeta(s+2)\zeta(-3-s)}{\zeta(-1-s)\zeta(s+4)}  \label{eqn70}
\end{equation}
This is a simple recursion relation for the $\zeta(s)$ with no other functions involved. It has five points ($s+1, s+2, s+3, s+4$ and $s-1$) around the point of interest $s$ and six points in other quadrants ($1-s, -1-s, -2-s,  -3-s, -s$ and $2-s$). These recursion relations are not symmetric in change of variable of $s\rightarrow1-s$ but have a symmetry of $s\rightarrow-1-s$ or $s\rightarrow1-s+c_1$ in general where $c_1\in{R}$.

\section{Functional Equations, Imaginary Increments}
By allowing an increment of $i\alpha$ in equation (\ref{eqn20}) yields the following ($\alpha \in R$)
\begin{equation}
\zeta(s)=\frac{s\zeta(1-s)\zeta(1+s)[i(s+i\alpha)tanh(\frac{\pi\alpha}{2})\zeta(1-s-i\alpha)\zeta(1+s+i\alpha)+\zeta(s+i\alpha)\zeta(-s-i\alpha)2\pi]}{2\pi\zeta(-s)[-2\pi{i}tanh(\frac{\pi\alpha}{2})\zeta(s+i\alpha)\zeta(-s-i\alpha)+(s+i\alpha)\zeta(1+s+i\alpha)\zeta(1-s-i\alpha)]}  \label{eqn80}
\end{equation}
This equation contains a parameter $\alpha$ which is formally set as  $\alpha \in R$. However, nothing is preventing it to behave as $\alpha \in C$. This is valid for all following equations involving the $\alpha$. 
By making an increment of $i\alpha$ in equation (\ref{eqn30}) yields the following.
\begin{equation}
0=16\pi^4\frac{\zeta(s+i\alpha)^2\zeta(-1-s-i\alpha)^2}{\zeta(1-s-i\alpha)^2\zeta(s+2+i\alpha)^2}+16\pi^4\frac{\zeta(s)^2\zeta(-1-s)^2}{\zeta(1-s)^2\zeta(s+2)^2} \ + \nonumber
\end{equation}
\begin{equation}
-32\pi^4\frac{\zeta(s)\zeta(-1-s)\zeta(s+i\alpha)\zeta(-1-s-i\alpha)}{\zeta(1-s)\zeta(s+2)\zeta(1-s-i\alpha)\zeta(s+2+i\alpha)}  \ + \nonumber 
\end{equation}
\begin{equation}
+8\pi^2\alpha^2\frac{\zeta(s+i\alpha)\zeta(-1-s-i\alpha)}{\zeta(1-s-i\alpha)\zeta(s+2+i\alpha)}+\alpha^4+\alpha^2  \ + \nonumber 
\end{equation}
\begin{equation}
+8\pi^2\alpha^2\frac{\zeta(s)\zeta(-1-s)}{\zeta(1-s)\zeta(s+2)} \label{eqn90}
\end{equation}
On the other hand, by managing the $s^2+s$ term while allowing the increment of $i\alpha$ in equation (\ref{eqn30}) leads to the following
\begin{equation}
\zeta(s+i\alpha)=\frac{\zeta(1-s-i\alpha)\zeta(2+s+i\alpha)[\zeta(s)\zeta(-1-s)-\frac{1}{4\pi^2}((2s+1)i\alpha-\alpha^2)\zeta(1-s)\zeta(s+2)]}{\zeta(-1-s-i\alpha)\zeta(1-s)\zeta(s+2)} \label{eqn100}
\end{equation}
We can eliminate $s$ in equation (\ref{eqn100}) and we will obtain the following.
\begin{equation}
\frac{-i\alpha}{2\pi^2}=\frac{\zeta(s+1+i\alpha)\zeta(-2-s-i\alpha)}{\zeta(-s-i\alpha)\zeta(s+3+i\alpha)}  \ + \nonumber
\end{equation}
\begin{equation}
-\frac{\zeta(s+1)\zeta(-2-s)}{\zeta(-s)\zeta(s+3)}  \ + \nonumber
\end{equation}
\begin{equation}
-\frac{\zeta(s+i\alpha)\zeta(-1-s-i\alpha)}{\zeta(1-s-i\alpha)\zeta(s+2+i\alpha)}  \ + \nonumber
\end{equation}
\begin{equation}
+\frac{\zeta(s)\zeta(-1-s)}{\zeta(1-s)\zeta(s+2)} \label{eqn110}
\end{equation}
Eliminating the left hand side in this equation yields the following.
\begin{equation}
0=\frac{\zeta(s+2+i\alpha)\zeta(-3-s-i\alpha)}{\zeta(-1-s-i\alpha)\zeta(s+4+i\alpha)}  \ + \nonumber
\end{equation}
\begin{equation}
-\frac{\zeta(s+2)\zeta(-3-s)}{\zeta(-1-s)\zeta(s+4)}  \ + \nonumber
\end{equation}
\begin{equation}
-2\frac{\zeta(s+1+i\alpha)\zeta(-2-s-i\alpha)}{\zeta(-s-i\alpha)\zeta(s+3+i\alpha)}  \ + \nonumber
\end{equation}
\begin{equation}
+2\frac{\zeta(s+1)\zeta(-2-s)}{\zeta(-s)\zeta(s+3)}  \ + \nonumber
\end{equation}
\begin{equation}
+\frac{\zeta(s+i\alpha)\zeta(-1-s-i\alpha)}{\zeta(1-s-i\alpha)\zeta(s+2+i\alpha)}  \ + \nonumber
\end{equation}
\begin{equation}
-\frac{\zeta(s)\zeta(-1-s)}{\zeta(1-s)\zeta(s+2)} \label{eqn120}
\end{equation}
This equation is the most general recursion relation for the $\zeta(s)$ over the complex plane with no other functions involved. It has nine points ($s+1, s+2, s+3, s+4, s+i\alpha, s+1+i\alpha, s+2+i\alpha, s+3+i\alpha, s+4+i\alpha$) around the point of interest $s$ and ten points in other quadrants ($1-s, -1-s,  -2-s, -s, -1-s-i\alpha, -2-s-i\alpha, 1-s-i\alpha, -3-s-i\alpha, -s-i\alpha$ and $-3-s$). If $\alpha\rightarrow-i$, this equation turns to equation (\ref{eqn70}). These recursion relations are not symmetric in change of variable of $s\rightarrow1-s$ but have various other symmetries depending on the complexity of the equation itself. In general the symmetry is of type $s\rightarrow1-s+c_2$ where $c_2\in{C}$.
\section{Additional Functional Equations and Results}
By starting from equation (\ref{eqn30}) we can iterate to yield
\begin{equation}
\zeta(s)=\frac{(-1)^{n+1}s(s+1)(s+2)...(s+2n)(s+2n+1)\zeta(1-s)\zeta(s+2n+2)}{(4{\pi}^2)^{n+1}\zeta(-s-2n-1)}  \label{eqn300}
\end{equation}
Here $n=0,1,2,3...$. To avoid apparent singularities, the following needs to be taken into account.
\begin{equation}
\zeta(-2N)=0, \ \ N=1,2,3...  \label{eqn302}
\end{equation}
and that the $\zeta(s)$ has a pole at $s=1$. This generates some conditions for the singularities. Some of them can be overcome by approaching through a proper limit. The index can grow to infinity. An interesting equation will follow if we set $s=\frac{1}{2}$.
\begin{equation}
\zeta(-\frac{3}{2}-2n)=\frac{(-1)^{n+1}(\frac{1}{2})(\frac{3}{2})(\frac{5}{2})(\frac{7}{2})...(\frac{1}{2}+2n)(\frac{1}{2}+2n+1)\zeta(\frac{5}{2}+2n)}{(4{\pi}^2)^{n+1}}  \label{eqn310}
\end{equation}
This offers an immediate link between the zeta function's values for certain positive and negative half-fractional arguments on the real axis.

One can start the iteration in a different way from (\ref{eqn30}) to obtain
\begin{equation}
\frac{\zeta(s+2)}{\zeta(-1-s)}=\frac{(-1)^{m+1}(4{\pi}^2)^{m+1}\zeta(s-2m)}{\zeta(2m+1-s)s(s+1)[(s-1)(s-2)...(s-2m+1)(s-2m)]} \label{eqn315}
\end{equation}
Here $m=0,1,2,3...$ and similar conditions as above need to be taken into account when singularities are to be avoided. The index can grow to infinity. If we set $s=\frac{-3}{2}$ we will get the same link (\ref{eqn310}) as above. 

By iterating equation (\ref{eqn20}) by its main functional factor, one will get 
\begin{equation}
\frac{\zeta(s)}{\zeta(1-s)}=\frac{s(s+1)(s+2)...(s+n)\zeta(n+1+s)\alpha{_n(s)}}{(2\pi)^{n+1}\zeta(-s-n)}  \label{eqn320}
\end{equation}
The function $\alpha{_n(s)}$ is a function of the iteration index $n$ and the argument $s$
\begin{equation}
\alpha{_n(s)}=\frac{1}{2}[i^{n}\cdot tan(\frac{\pi{s}}{2})((-1)^{n+1}-1)+i^{n+3}\cdot ((-1)^{n}-1)]    \label{eqn330}
\end{equation}
Setting $s=\frac{1}{2}$ in (\ref{eqn320}) yields
\begin{equation}
\zeta(-\frac{1}{2}-n)=\frac{(\frac{1}{2})(\frac{3}{2})(\frac{5}{2})(\frac{7}{2})...(\frac{1+2n}{2})\zeta(\frac{3}{2}+n)\alpha{_n(\frac{1}{2})}}{(2\pi)^{n+1}}  \label{eqn335}
\end{equation}
and
\begin{equation}
\alpha{_n(\frac{1}{2})}=\frac{1}{2}[i^{n}\cdot ((-1)^{n+1}-1)+i^{n+3}\cdot ((-1)^{n}-1)]    \label{eqn340}
\end{equation}
This is another link between the zeta function's values for positive and negative half-fractional arguments on the real axis. It is complementary to the equation (\ref{eqn310}) being twice as dense. More iterated recursion relations can be developed along these lines. 

As examples of the equations above consider the following. By setting $s=-\frac{1}{2}$ in equation (\ref{eqn30}) one can obtain
\begin{equation}
\frac{\zeta(-\frac{1}{2})}{\zeta(\frac{3}{2})}=\pm\frac{1}{4{\pi}}  \label{eqn380}
\end{equation}
The negative sign is valid. From equation (\ref{eqn310}) one would get after setting $n=0$ and $n=1$ 
\begin{equation}
\zeta(-\frac{3}{2})=\frac{(\frac{-1}{2})(\frac{3}{2})\zeta(\frac{5}{2})}{4{\pi}^2}  \label{eqn390}
\end{equation}
\begin{equation}
\zeta(-\frac{7}{2})=\frac{(\frac{1}{2})(\frac{3}{2})(\frac{5}{2})(\frac{7}{2})\zeta(\frac{9}{2})}{(4{\pi}^2)^2}  \label{eqn400}
\end{equation}
From equation (\ref{eqn335}) for $n=2$ 
\begin{equation}
\zeta(-\frac{5}{2})=\frac{(\frac{1}{2})(\frac{3}{2})(\frac{5}{2})\zeta(\frac{7}{2})}{(2{\pi})^3}  \label{eqn500}
\end{equation}

\section{Discussion}
All of the recursion relations obtained in this work are equivalent to the original functional equation (\ref{eqn10}). This can be verified by placing the equation (\ref{eqn10}) to any of them and observing it to become an identity. Thus each is a functional equation for the Riemann zeta function. The new recursion relations (\ref{eqn30}...\ref{eqn300}, \ref{eqn315} and \ref{eqn320}) offer links between points of the zeta function in the complex plane or along the real axis (\ref{eqn310} and \ref{eqn335}). They are useful for evaluating the Riemann's zeta function at a point when values at some other points (dictated by the recursion) are known. These recursion relations are not symmetric in change of variable of $s\rightarrow1-s$ but have various other symmetries depending on complexity of the relation itself. In general the symmetry is of type $s\rightarrow1-s+c_2$ where $c_2\in{C}$.

The equation (\ref{eqn30}) has a similar symmetric structure as the original functional equation (\ref{eqn10}), namely 
\begin{equation}
\frac{\zeta(s)}{\zeta(1-s)}=-s(s+1)\frac{\zeta(s+2)}{4\pi^2\zeta(-1-s)}  \label{eqn600}
\end{equation}
The combination on the left is common to many of the expressions above with varying arguments. It is here also on the right side. Another common combination in these functional equations is the left hand side of
\begin{equation}
\frac{\zeta(s)\zeta(-1-s)}{\zeta(1-s)\zeta(s+2)}=-s(s+1)\frac{1}{4\pi^2}  \label{eqn610}
\end{equation}
being a simple analytic function.

There seems to be a great number of recursion relations for the zeta function, possibly an infinity. The iterated functional equations point to the fact that there may exist families of recursion relations. Often the simplest ones merge to one of the equations in the preceding chapters by a simple transformation. Several new relations can be derived from equation (\ref{eqn10}) or equation (\ref{eqn30}) in different ways or by combining any of the relations presented. The author considers equation (\ref{eqn30}) as the central result of this work.

\end{document}